\documentclass[11pt]{article}

\usepackage{latexsym,geometry}
\usepackage{amsmath,amsfonts,amssymb}
\usepackage{graphicx,color,epsfig}

\def\R{\mathbb{R}}
\def\N{\mathbb{N}}

\def\leq{\leqslant}
\def\geq{\geqslant}

\newtheorem{Def}{Definition}[section]
\newtheorem{Thm}[Def]{Theorem}
\newtheorem{Lem}[Def]{Lemma}
\newtheorem{Pro}[Def]{Proposition}
\newtheorem{Cor}[Def]{Corollary}
\newtheorem{Claim}[Def]{Claim}

\begin{document}

\title{Population Dynamics of Globally Coupled Degrade-and-Fire Oscillators}
\author{Alex Blumenthal$^1$ and Bastien Fernandez$^2$\footnote{On leave from Centre de Physique Th\'{e}orique, CNRS - Aix-Marseille Universit\'e - Universit\'e de Toulon, Campus de Luminy, 13288 Marseille CEDEX 09 France}}
\maketitle

\begin{center}
$^1$ Courant Institute of Mathematical Sciences\\
New York University\\
New York, NY 10012, USA

\vskip0.5cm

$^2$ Laboratoire de Probabilit\'es et Mod\`eles Al\'eatoires\\
CNRS - Universit\'e Paris 7 Denis Diderot\\
75205 Paris CEDEX 13 France\\
\end{center}

%
%
%
%
%
%
%
%
%
%
%

\begin{abstract}
This paper reports the analysis of the dynamics of a model of pulse-coupled oscillators with global inhibitory coupling. The model is inspired by experiments on colonies of bacteria-embedded synthetic genetic circuits. The total population can be either of finite (arbitrary) size or infinite, and is represented by a one-dimensional profile. Profiles can be discontinuous, possibly with infinitely many jumps. Their time evolution is governed by a singular differential equation. We address the corresponding initial value problem and characterize the dynamics' main features. In particular, we prove that trajectory behaviors are asymptotically periodic, with period only depending on the profile (and on the model parameters). A criterion is obtained for the existence of the corresponding periodic orbits, which reveals the existence of a sharp transition as the coupling parameter is increased. The transition separates a regime where any profile can be obtained in the limit of large times, to a situation where only trajectories with sufficiently large groups of synchronized oscillators perdure.  
\end{abstract}




\section{Introduction}
To determine the amount of collective order and how its depends on parameters is a central question in the analysis of systems of coupled oscillators \cite{PRK01}. A typical example is the Kuramoto model \cite{Kur84} where numerics have indicated that, while the oscillators evolve independently at weak coupling, as soon as the interaction strength exceeds a threshold, the overall collective behavior becomes more and more coherent when the coupling increases \cite{ABP-VRS05}. This phenomenology has been identified since the late 1970's. Yet, its full mathematical proof still remains to be achieved and only preliminary results have been obtained so far (see \cite{Strogatz00} for a summary of results and technical challenges).  

In the Kuramoto model, trajectories are smooth and most of the difficulties come from the presence of heterogeneities ({\sl i.e.}\ the individual oscillators' frequencies are randomly drawn), not to mention nonlinearities. However, in other circumstances, modelling results in systems with temporal singularities, namely bursts or spikes. This is typically the case of pulse-coupled oscillators such as the integrate-and-fire model \cite{P75}. 

In assemblies of pulse-coupled oscillators with excitatory couplings, global synchrony has been proved to hold for any interaction strength, not only in the homogeneous case \cite{MS90}, but also for certain heterogeneous models with distributed individual frequencies, thresholds and/or coupling parameters \cite{SU00}. For inhibitory couplings, full synchrony usually fails and instead, several distinct groups of synchronous units are observed (see however \cite{MS13} for a similar phenomenology in the case of excitatory couplings). However, mathematical proofs are scarce in this setting, except for populations consisting of two units \cite{EPG95}.

The current paper reports the mathematical analysis of a model of pulse-coupled oscillators with inhibitory coupling, inspired from a series of experiments on synthetic genetic circuits embedded in E.\ coli \cite{DM-PTH10}. The model mimics {\sl in vitro} dynamics of assemblies of interacting self-repressor genes. The analysis here not only addresses populations of arbitrary finite size throughout the coupling parameter range, but also continua of infinite populations.

As phenomenology is concerned, the critical feature of this model is the existence a sharp transition when the coupling strength increases. Below the threshold, every possible population distribution can be observed at large time, upon the choice of initial condition. Beyond that point, only distributions of grouped oscillators can persist - the stronger the interaction, the larger the groups. Any initially isolated oscillator must eventually join, and must remain, in a group. 

These features have been identified and mathematically justified for finite size populations in \cite{FT11,FT13}, except for the asymptotic periodicity of every trajectory (Theorem \ref{GLOBCONV} below), whose proof in \cite{FT11} turns out to be incomplete. Here, we provide a complete proof of this statement and we address the initial value problem. Following \cite{AvV93,B99,K91,MS13}, we also extend the analysis to the continuum of oscillators. In practice, this means dealing with a (singular) differential equation whose variable is a discontinuous real function with possibly infinitely many jumps (and consists of finitely many plateaus in the case of finite size populations). We also consider the initial value problem in the infinite-dimensional context and characterize the asymptotic dynamics, at least in the weak coupling regime. 

The paper is organised as follows. After the definition of the model (section \ref{S-DEF}), we describe in section \ref{S-BASIC} the basic features of solutions, mainly the fact that repeated firings must occur in every cell. Thanks to the mean-field nature of the coupling, the dynamics commutes with permutations of individuals in the population. We analyze the consequences of this symmetry on collective behaviors and grouping properties. In section \ref{S-INITVAL}, we consider the initial value problem and prove existence and uniqueness of global solutions for every piecewise constant initial condition, and also for every initial configuration when the coupling is sufficiently small. 
Section \ref{S-ASYMP} focuses on asymptotic properties of the dynamics. We first characterize the simplest periodic orbits and establish a necessary and sufficient condition on the coupling parameter for their existence. The analysis of this condition shows that an abrupt transition occurs as this parameter crosses a threshold. The transition separates a regime where all possible periodic trajectories exist, to a situation where only trajectories with sufficiently large plateaus perdure. In the same section, we prove that every finite-dimension solution must be asymptotically periodic, and usually approaches one of the previously mentioned periodic trajectories. We also prove that the same results hold for every solution, provided again that the coupling is not too large. Finally, we conclude the paper by a series of comments on robustness with respect to changes in the model. 

\section{Definitions}\label{S-DEF}
For simplicity, we assume that a single gene is involved in our process (instead of two genes in the original experiment \cite{DM-PTH10}) and that intercellular coupling is co-repressive and of mean-field type (as opposed to co-excitatory and spatially localized as in the original modelling \cite{MBHT09}). 

Each oscillator represents a self-repressor gene \cite{TD-A90} embedded in a host cell (same gene in every cell) labelled by a real number $x\in (0,1]$. The oscillator state at time $t\in\R^+$ is characterised by a real number $u(x,t)\in [0,1]$ that represents the so-called gene expression level (normalized concentration) in cell $x$ \cite{TCN01}. 

Intercellular coupling materializes via a repressor field. As in the original experiment, we consider that gene products involved in the feedback loop are small enough so that they can diffuse through membranes. Assuming in addition that strong stirring holds in the container, we also consider that mixing of diffused gene levels occurs very rapidly in the medium ({\sl i.e.}\ on a much shorter time scale than protein degradation processes). Accordingly, the repressor level $Mu(x,t)$ in cell $x$, at time $t$, is given by the following linear combination
\[
M u(x,t) = (1-\epsilon\eta)u(x,t)+\epsilon\eta\int_0^1u(y,t)dy,
\]
where ${\displaystyle\int_0^1}u(y,t)dy$ represents the mean field in the intercellular medium and where the diffusion coefficient $\epsilon\eta$ is composed by the product of the {\bf coupling parameter} $\epsilon$ with the {\bf threshold parameter} $\eta$ (the choice of $\epsilon\eta$ instead of only $\epsilon$ turns out to be more convenient in our results). In principle, $\eta$ can take any value in $(0,1)$, but we assume its value is (very) small, in agreement with the original observations. The parameter $\epsilon$ is restricted in a way that repressor levels never become negative, {\sl i.e.}\ we impose $0< \epsilon<1/\eta$.\footnote{The dynamics can also be defined for $\epsilon=0$ and is obvious in this case. The assumption $\epsilon>0$ is more convenient for the analysis.}   

With these preliminary definitions provided, the evolution rule can be given. The dynamics of gene expression levels is governed by the following singular differential equation\footnote{We use the notation $u(x,t-0):=\lim\limits_{s\to t,s<t}u(x,s)$. A similar definition holds for $u(x,t+0)$. Moreover, the sign symbol $\text{Sgn}$ is defined on $\R^+$ by
\[
\text{Sgn}(u)=1\ \text{if}\ u>0\quad \text{and}\quad \text{Sgn}(0)=0.
\]
}
\begin{equation}
\begin{aligned}
&\begin{aligned}
&\partial_t u(x,t)=-\text{Sgn}(u(x,t))&\text{if}&&Mu(x,t)>\eta\\ 
&\left\{\begin{array}{l}
u(x,t)=u(x,t-0)\\
u(x,t+0)=1
\end{array}\right.&\text{if}&&Mu(x,t)\leq\eta
\end{aligned}\quad \forall x\in (0,1],t\in (0,+\infty)\\
&\text{and}\ u(x,0)=u(x)\quad \forall x\in (0,1].
\end{aligned}
\label{DEFDYNAM}
\end{equation}
Throughout the paper, functions depending on a single variable depend on $x$ (unless otherwise stated) and are viewed as one-dimensional {\bf profiles}. As a consequence, no confusion results from using the symbol $u(x)$ to denote the initial condition of the solution $u(x,t)$.

Equation \eqref{DEFDYNAM} is inspired by the delay-differential equation that has been introduced in \cite{MBHT09} in order to reproduce the experimental oscillations. Both our model and the one in \cite{MBHT09} obey similar principles: 
\begin{itemize}
\item[$\bullet$] A significant repressor level ($M u(x,t)>\eta$) in a cell prevents any production. In this case, the corresponding gene expression level decays slowly (constant speed $-1$) due to degradation. 
\item[$\bullet$] When the repressor expression becomes (very) small ($M u(x,t)\leq\eta$), repression can no longer prevent production and the gene is synthetized at fast rate (also on infinitesimal scales when compared to degradation processes) until saturation is reached. This instantaneous production event is called a {\bf firing}. (Lemma \ref{BOUNDM} below actually shows that we have $M u(x,t)\geq\eta$ for all $(x,t)$; hence firings occur exactly when $M u(x,t)=\eta$.) 
\item[$\bullet$] Gene levels may locally reach and stay at negligible values ($u(x,t)=0$) while their repressor level remains above threshold. In this case, the next firing is delayed and cell $x$ will evolve in sync with any other cell $y$ whose expression level has reached zero in the mean time (a property called 'grouping process', see section \ref{S-GROUPROCESS} below).
\end{itemize}

The assumptions on the initial {\bf profile} $u$ are as follows
\begin{itemize}
\item[(1)] $u$ is a real Borel measurable function (with values in $(0,1]$); hence the quantity $Mu(x):=Mu(x,0)$ is well-defined for every $x\in (0,1]$
\item[(2)] $u$ is non-decreasing on $(0,1]$. No loss of generality results from making this assumption because the dynamics commutes with label exchange;\footnote{Namely, if the relation 
\[
v(x,t)=u(x,t) ,\ \forall x\neq x_1,x_2 \quad  \text{and}\quad \left\{\begin{array}{l}
v(x_1,t)=u(x_2,t)\\
v(x_2,t)=u(x_1,t)
\end{array}\right.
\]
holds for $t=0$, and if $u(x,t)$ is a solution of the equation \eqref{DEFDYNAM}, then $v(x,t)$ is also a solution and the previous relation holds for all $t\in (0,+\infty)$.} hence the cell can always be labelled in a way that their expression level are monotonically ordered.
\item[(3)] $u$ be {\bf c\`agl\`ad}, {\sl i.e.}\ $u$ is also left continuous on $(0,1]$, the  existence of a right limit at every point of $(0,1)$ follows from monotonicity.
\item[(4)] $Mu(0+0)>\eta$ and $u(0+0)<u(1)= 1$.
\end{itemize}
The following comments are in order
\begin{itemize}
\item[$\bullet$] the condition $Mu(0+0)>\eta$ ensures the existence and uniqueness of the solution $u(x,t)$ locally for $t>0$ in a neighborhood of 0 (notice again that we will show that any solution satisfies $M u(x,t)\geq\eta$ for all $(x,t)$).
\item[$\bullet$] the condition $u$ is c\`agl\`ad is consistent with the same feature of the solution profile at every time (see expression \eqref{EXPRSOLUTION} and Proposition \ref{PROPFIR} below)
\item[$\bullet$] the condition $u(1)=1$ states that the cell(s) with highest expression level in the initial population has (have) just fired at $t=0$. Properties of the firing times in Lemma \ref{PROPFIR} below imply that every solution satisfies this property at some moment in time (at infinitely many moments indeed); hence that assumption amounts to a time translation.
\item[$\bullet$] the inequality $u(0+0)<u(1)$ is also a matter of convenience. It ensures that the initial population is not in full synchrony; otherwise the dynamics would be trivial (single oscillator) and does not require any elaborated investigation. 
\end{itemize}
We shall require more conditions below when we address the existence of global solutions ({\sl i.e.}\ functions $u(x,t)$ that satisfy \eqref{DEFDYNAM} for all $(x,t)\in (0,1]\times \R^+$, for given $\epsilon$ and $\eta$). 

\section{Basic dynamical features}\label{S-BASIC}
Postponing the existence of global solutions of equation \eqref{DEFDYNAM} to section \ref{S-INITVAL} below, in this section, we describe basic but essential solution features. On these properties depend both the analysis and the formulation of the results in the sections below. We begin with temporal features of individual oscillators before passing to the description of the collective properties of the population dynamics.

\subsection{Properties of the firing times}
Here, focus is made on basic temporal features; in particular on the facts that (see Lemma \ref{PROPFIR})
\begin{itemize}
\item[$\bullet$] every cell must fire repeatedly forever, and 
\item[$\bullet$] between every two consecutive firings in $x$, there must be exactly one firing in every other cell $y\neq x$, unless $y$ fires simultaneously with $x$.
\end{itemize}
A first statement justifies the firing time definitions below.
\begin{Lem}
For every solution $(x,t)\mapsto u(x,t)$, the inequality $Mu(x,t)\geq \eta$ holds for all $(x,t)\in (0,1]\times \R^+$. 
\label{BOUNDM}
\end{Lem}
{\sl Proof.}\ The proof is by contradiction. It relies on the fact that, for every $x\in (0,1]$, the function $t\mapsto u(x,t)$ is c\`agl\`ad. 
This property follows trivially from equation \eqref{DEFDYNAM} when $M u(x,t)\leq\eta$ and is a consequence of the fact that $t\mapsto u(x,t)$ must be continuous when $M u(x,t)>\eta$. 
Together with Lebesgue's dominated convergence theorem, it implies that each function $t\mapsto M u(x,t)$ is also c\`agl\`ad.

\noindent
By contradiction, assume the existence of $(x_0,t_0)$ such that $Mu(x_0,t_0)<\eta$. We must have $t_0>0$ due to the initial assumption $M u(x,0)>\eta$ for all $x\in (0,1]$. Then left continuity implies the existence of $t_1<t_0$ such that $Mu(x_0,t)<\eta$ and hence $u(x_0,t)=1$ for all $t\in (t_1,t_0]$. Since every expression level is smaller or equal to 1, it follows that $Mu(x,t)<\eta$ and hence $u(x,t)=1$ for all $(x,t)\in (0,1]\times (t_1,t_0]$. Using the definition of $M$, the latter yields $Mu(x,t)=1$ for all $(x,t)\in (0,1]\times (t_1,t_0]$, which contradicts the assumption $Mu(x_0,t_0)<\eta$. 
\hfill
\hfill $\Box$

Lemma \ref{BOUNDM} implies that firing events occur exactly when the repressor level reaches the threshold $\eta$. In particular, the {\bf first firing time} in cell $x$, defined as 
\[
T_1u(x):=\sup\{t> 0\ :\ M u(x,s)>\eta,\ \forall 0\leq s<t\},
\]
can be characterized as follows
\[
T_1u(x)=\inf\{t> 0\ :\ M u(x,t)=\eta\}.
\]
The c\`agl\`ad property of $t\mapsto M u(x,t)$ and the initial assumption $M u(0+0,0)>\eta$ ensure that $T_1u(x)>0$ for every $x\in (0,1]$. Anticipating that $T_1u(x)<+\infty$ and that the repressor level immediately after firing, {\sl i.e.}\ $Mu(x,T_1u(x)+0)>\eta$, lies above $\eta$, the second firing time can be defined similarly, namely
\[
T_2u(x)=\inf\{t> T_1u(x)\ :\ M u(x,t)=\eta\}.
\]
Repeating the argument, one defines the successive firing times as follows 
\[
T_{n+1}u(x)=\inf\{t>T_nu(x)\ :\ Mu(x,t)=\eta\},\ \forall n\in\N,
\]
and anticipating also on the fact that $T_nu(x)\xrightarrow{n\to+\infty} +\infty$, one obtains the following explicit expression of global solutions\footnote{$u^+:=\max\{u,0\}$.}
\begin{equation}
u(x,t)=\left\{\begin{array}{ccl}
(u(x)-t)^+&\text{if}&0\leq t\leq T_1u(x)\\
(1-t+T_nu(x))^+&\text{if}&T_nu(x)< t\leq T_{n+1}u(x),\ n\in\N
\end{array}\right.\ \forall x\in (0,1].
\label{EXPRSOLUTION}
\end{equation}
All required assumptions above are listed in the main statement of this section, which we now formulate.
\begin{Pro}
For every solution $(x,t)\mapsto u(x,t)$, the following properties hold.
\begin{itemize}
\item[$\bullet$] For every $x\in (0,1]$ and $n\in\N$, the firing time $T_nu(x)$ is well-defined and is finite.
\item[$\bullet$] Every function $x\mapsto T_nu(x)$ is non-decreasing and left continuous (and therefore c\`agl\`ad) and we have $T_nu(1)\leq T_{n+1}u(0+0)$.
\item[$\bullet$] For every $x\in (0,1]$ we have $T_1u(x)\geq Mu(x)-\eta$ and $T_{n+1}u(x)-T_nu(x)\geq (1-\epsilon\eta)(1-\eta)$ for all $n\in\N$.
\end{itemize}
\label{PROPFIR}
\end{Pro}
{\sl Proof.}\ 
Most of the effort consists in proving the statement for $n=1$ (and that $T_2$ is well-defined; we already know that $T_1$ is well-defined); the other cases will follow by induction, by applying these conclusions to the solution at appropriate successive times. 

\noindent
As we shall see, the proof for $n=1$ only relies on the following assumptions on $u$
\begin{equation}
u\ \text{is non-decreasing and left continuous},\ u(1)= 1,\ \text{and}\ Mu(x)>\eta,\ \forall x\in (0,1],
\label{ASSUMPROFILE}
\end{equation}
and does not need that $u(0+0)<u(1)$, {\sl i.e.}\ that the cells are initially out of sync. 
\smallskip

\noindent
$\bullet$ {\sl Proof of monotonicity of the function $T_1$}. 
Given two arbitrary points $x_1<x_2$, using that $\min\{T_1u(x_1),T_1u(x_2)\}\in \R^+_\ast\cup\{+\infty\}$, 
let $t>0$ be 
such that $t\leq \min\{T_1u(x_1),T_1u(x_2)\}$. Expression \eqref{EXPRSOLUTION}, together with $u(x_1)\leq u(x_2)$, implies $u(x_1,t)\leq u(x_2,t)$ and then $Mu(x_1,t)\leq Mu(x_2,t)$ from where the inequality $T_1u(x_1)\leq T_1u(x_2)$ follows. 

\noindent
$\bullet$ {\sl Proof of the inequality $T_1u(1)\leq T_2u(0+0)$}. 
The definition $u(x,T_1u(x)+0)=1$ implies that 
\[
Mu(x,T_1u(x)+0)\geq Mu(y,T_1u(x)+0),\ \forall y\neq x.
\]
Using also monotonicity of $T_1$, this implies that the second firing $T_2u(x)$ in cell $x$ cannot happen before cell 1 has first fired, {\sl i.e.}\ $T_1u(1)\leq T_2u(x)$ for all $x\in (0,1]$. The inequality $T_1u(1)\leq T_2u(0+0)$ immediately follows. (Of note, if $T_1u(x)=\infty$ for some $x$, then there is nothing to prove. Moreover, the same argument, together with monotonicity of $T_1$, implies monotonicity of the function $T_2$.) 
\smallskip

For the proof of finiteness of $T_1$, we shall rely on the following statement. 
\begin{Lem} {\rm (i)} $T_1u(0+0)<1$. 

\noindent
{\rm (ii)} Assume that $T_1u(x_1)<+\infty$ for some $x_1\in (0,1]$. Then, there exists $x_2\neq x_1$ and $n\in\{1,2\}$ such that $T_1u(x_1)<T_nu(x_2)<T_1u(x_1)+1$. 
\label{BASICCLAIM}
\end{Lem}
{\sl Proof of the Lemma.} (i) By contradiction, assume that $T_1u(0+0)\geq 1$. Then monotonicity of $T_1$ implies $T_1u(x)\geq 1$ for all $x$. Expression \eqref{EXPRSOLUTION}, together with $u(x)\leq 1$, yields $u(x,1)=0$ for all $x$. This in turns gives $Mu(x,1)=0$, which is impossible by Lemma \ref{BOUNDM}.

\noindent
(ii) The arguments are similar. By contradiction, given that $T_1u(x_1)<+\infty$, assume that we have 
\[
T_1u(x)\geq T_1u(x_1)+1, \forall x\in (x_1,1]\quad\text{and}\quad T_2u(x)\geq T_1u(x_1)+1, \forall x\in (0,x_1].
\]
Then, as for statement (i), we conclude that $Mu(x,T_1u(x_1)+1)=0$ for all $x$, which is impossible. 
Statement (ii) is nothing but a condensed formulation of the assumption negation.
\hfill $\Box$
\smallskip

\noindent
$\bullet$ {\sl Proof of finiteness of the function $T_1$}. By contradiction, assume that $T_1u(x)=+\infty$ for some $x\in (0,1]$. Then, monotonicity of $T_1$ and Lemma \ref{BASICCLAIM} imply the existence of $x_0\in (0,1]$ such that
\[
\left\{\begin{array}{ccl}
T_1u(x)<+\infty&\text{if}&0<x<x_0\\
T_1u(x)=+\infty&\text{if}&x_0\leq x\leq 1
\end{array}\right.
\]
We consider two cases; either $T_1u(x_0-0)<+\infty$ or $T_1u(x_0-0)=+\infty$. In the first case, the property $T_1u(1)\leq T_2u(0+0)$ and monotonicity of $T_2u$ imply that no firing can occur in a cell after time $T_1u(x_0-0)$. Using similar arguments as in the proof of Lemma \ref{BASICCLAIM}, this implies $Mu(x,T_1u(x_0-0)+1)=0$ for all $x$, which is impossible. 

\noindent
Assume now that $T_1u(x_0-0)=+\infty$. After time $T_1u(x_0-\frac{1}{2\epsilon})$, only those cells in the interval $(x_0-\frac{1}{2\epsilon},x_0)$ can fire. Therefore, we have 
\[
\left\{\begin{array}{ccl}
u(x,T_1u(x_0-\frac{1}{2\epsilon})+1)=0&\text{if}&x\not\in (x_0-\frac{1}{2\epsilon},x_0)\\
u(x,T_1u(x_0-\frac{1}{2\epsilon})+1)\leq 1&\text{if}&x\in (x_0-\frac{1}{2\epsilon},x_0)
\end{array}\right.
\]
which implies  
\[
Mu(x,T_1u(x_0-\frac{1}{2\epsilon})+1)\leq \epsilon\eta\frac{1}{2\epsilon}<\eta, \forall x\not\in (x_0-\frac{1}{2\epsilon},x_0)
\]
which is again impossible. 
\smallskip

\noindent
$\bullet$ {\sl Proof of left continuity of the function $T_1$}. 
Monotonicity obviously implies that the limit $T_1u(x-0)$ exists,  and $T_1u(x-0)\leq T_1u(x)$, for every $x\in (0,1]$. Fix $x\in (0,1]$ and choose any $y\in (0,x)$. We have 
\[
M u(x,T_1u(y))= M u(y,T_1u(y))+(1-\epsilon\eta)\left(\left(u(x)-T_1u(y)\right)^+-\left(u(y)-T_1u(y)\right)^+\right)
\]
Using $M u(y,T_1u(y))=\eta$ and left continuity of the initial profile $u$, we obtain
\[
\lim_{y\to x^-}M u(x,T_1u(y))=\eta.
\]
Furthermore, that $t\mapsto Mu(x,t)$ is c\`agl\`ad for every $x\in (0,1]$ (see proof of Lemma \ref{BOUNDM}) implies that the limit here is equal to $Mu(x,T_1u(x-0))$. Hence, $Mu(x,T_1u(x-0))=\eta$ and the first firing time definition implies $T_1u(x-0)\geq T_1u(x)$ from where the desired conclusion $T_1u(x-0)= T_1u(x)$ follows. 
\smallskip

\noindent
$\bullet$ {\sl Proof of the inequality $T_{2}u(x)-T_1u(x)\geq (1-\epsilon\eta)(1-\eta)$ for all $x\in (0,1]$}. (The inequality $T_1u(x)\geq Mu(x)-\eta$ is a direct consequence of the initial assumption $Mu(x)>\eta$ together with $Mu(x,t)\geq Mu(x)-t$.) First, notice that the expression level in a cell that is about to fire, namely $u(x,T_1u(x))$, must minimise the expression levels at this time 
and must not lie above $\eta$. 
We conclude that
\[
Mu(x,T_1u(x)+0)\geq Mu(x,T_1u(x))+(1-\epsilon\eta)(1-u(x,T_1u(x)))\geq \eta+(1-\epsilon\eta)(1-\eta),
\]
from where the desired inequality immediately follows. As a by-product, this inequality also implies that $Mu(x,T_1u(x)+0)>\eta$ for all $x\in (0,1]$; hence the second firing time function $T_2$ is indeed well-defined, as claimed.
\smallskip

\noindent
$\bullet$ {\sl Induction step}. At this stage, it remains to show that (an appropriate translation of) the population profile immediately after cell 1 has fired satisfies the same assumptions \eqref{ASSUMPROFILE} as the initial function $u$. 
There are two cases; either $T_1u(1)<T_2u(x)$ for all $x\in (0,1]$ or $T_1u(1)=T_2u(x_0)$ for some $x_0\in (0,1]$. In the first case, we consider the limit function $u_1(x):=u(x,T_1u(1)+0)$ for all $x$. The assumption $T_1u(1)<T_2u(x)$, together with expression \eqref{EXPRSOLUTION}, implies (notice that all expression levels must be positive immediately after firing)
\[
u_1(x)=1-T_1u(1)+T_1u(x),\ \forall x\in (0,1],
\]
and so monotonicity and left continuity of $T_1$ imply the same properties for $u_1$ (and we obviously have $u_1(1)= 1$). Moreover, the same assumption also implies $Mu(x,T_1u(1))>\eta$, and {\sl a fortiori} $Mu(x,T_1u(1)+0)>\eta$, for all $x$ such that $T_1u(x)<T_1u(1)$.\footnote{Indeed, we have $T_2u(x)>T_1u(x)$ for all $x\in (0,1]$.} Therefore, the function $u_1$ satisfies all the assumptions of \eqref{ASSUMPROFILE} and the induction can proceed. 
\noindent

In the case where $T_1u(1)=T_2u(x_0)$ for some $x_0\in (0,1)$, we need to apply some spatial translation to the profile after the firing at time $T_1u(1)$ in order to obtain a monotonic function on $(0,1]$. Towards that goal, let 
\[
x_\text{max}:=\sup\{x\in (0,1)\ :\ T_2u(x)=T_1u(1)\}.
\]
Notice that the proof of Lemma \ref{BASICCLAIM} can be repeated for the function $T_2$ to conclude that $T_2(x)<+\infty$ for all $x$. Then the proof of left continuity of $T_1$ applies {\sl mutatis mutandis} to prove that the function $T_2$ must also be left continuous. Using also $T_2u(1)\geq T_1u(1)+(1-\epsilon\eta)(1-\eta)$, it follows that $x_\text{max}<1$ and then $T_2(x_\text{max})=T_1u(1)$, and we set
\[
u_1(x)=\left\{\begin{array}{cclcl}
u(x+x_\text{max},T_1u(1)+0)&=&1-Tu_1(1)+T_1u(x+x_\text{max})&\text{if}&0<x\leq 1-x_\text{max}\\
u(x+x_\text{max}-1,T_1u(1)+0)&=&1&\text{if}&1-x_\text{max}<x\leq 1
\end{array}\right.
\]
This function $u_1$ is obviously non-decreasing, left continuous and we have $u_1(1)=1$. Moreover, the definition of $x_\text{max}$ implies $T_2u(x)>T_1u(1)$ for all $x\in (x_\text{max},1]$; hence we have
\[
Mu(x,T_1u(1)+0)=Mu(x,T_1u(1))>\eta,\ \forall x\in (x_\text{max},1]\ :\ T_1u(x)<T_1u(1).
\]
As before, this easily implies $Mu_1(x)>\eta$ for all $x\in (0,1]$ and all assumptions of \eqref{ASSUMPROFILE} hold. The induction can proceed and the proof of Proposition \ref{PROPFIR} is complete. 
\hfill
\hfill $\Box$

\subsection{Grouping properties}\label{S-CLUSTER}
In this section we review grouping properties of the dynamics. These properties are consequences of the label exchange symmetry. For simplicity, the properties are formulated in terms of 
$T_1u$. As in the proof of Proposition \ref{PROPFIR} above, they extend to every finite time, by induction on profiles, after every full cycle of firings.

\subsubsection{Group invariance} If $u(x)=u(y)$ then $T_1u(x)=T_1u(y)$. 

\noindent
Therefore, if at some time $t_1$, $u(x,t_1)$ is constant on some interval, then it remains constant on this interval for all $t>t_1$. In other words, cells holding the same gene expression level define a {\bf group} \cite{MS90} and evolve in unison. ("Cluster" is another term for such groups \cite{ABP-VRS05}.) 

\subsubsection{Firing without grouping}\label{ONE} If $T_1u(x)\leq u(x)$, then $T_1u(x)<T_1u(y)$ for all $y$ such that $u(x)<u(y)$. 

\noindent
Therefore, if a cell $x$ (or a group of cells including $x$) fires before its level has hit 0, then it does so unaccompanied by any cell (or group) whose expression level differs from $x$ at this instant.

\subsubsection{No grouping regime}\label{TWOPRIME} If $\epsilon\leq 1$, then for every initial profile $u$, we have $T_1u(x)\leq u(x)$ for all $x\in (0,1]$. As a consequence, no grouping can occur for $\epsilon \in (0,1]$. 

\noindent
We prove this property by contradiction. Assume that $\epsilon\leq 1$ and $T_1u(x)>u(x)$ for some $x\in (0,1]$. 
Then, using the expression of the solution prior to the first firing (see equation \eqref{EXPRSOLUTION}) and the inequality $\int_0^1 u(y,t)dy\leq 1$, we get 
\[
Mu(x,u(x))=\epsilon\eta\int_0^1 u(y,u(x))dy\leq \eta\ \text{when}\ \epsilon\leq 1,
\]
which, considering that $Mu(x,0)>\eta$, is incompatible with $T_1u(x)>u(x)$.

\subsubsection{Grouping process}\label{S-GROUPROCESS} If  $u(x)<u(y)<T_1u(x)$, then $T_1u(y)=T_1u(x)$. 

\noindent
If the expression level in a cell/group hits 0 before firing, then the cell/group joins any other cell/group whose level also hits 0 before the same firing. 

\subsubsection{Maximal size of a forming/inflating group}\label{THREEPRIME} When $\epsilon>1$, the maximal size of a forming/inflating group before a firing is $1-\tfrac{1}{\epsilon}$.

\noindent
To see this, by contradiction again, assume there exists $t$ such that $u(x,t)=0$ for all $x\in (0,y]$ with $y>1-\tfrac{1}{\epsilon}$. Then, using that $u\leq 1$ for the rest of cells, we would have 
\[
Mu(x,t)\leq \epsilon\eta(1-y)<\epsilon\eta\frac{1}{\epsilon}=\eta,
\]
which is impossible.

Group invariance and the grouping process imply that the total plateaus' length ({\sl i.e.}\ the Lebesgue measure of the set where $u(\cdot,t)$ is constant) is a non-decreasing function of time. Since this length cannot exceed 1, it must converge. In other words, the total length of intervals where grouping occurs upon firings must vanish as $t\to +\infty$. 

A by-product of the properties in section \ref{TWOPRIME} and \ref{THREEPRIME} above is that full grouping of cells ({\sl i.e.}\ complete synchrony) can never be achieved before any firing (and hence in finite time), unless all cells are initially in sync. 

\section{Analysis of the initial value problem}\label{S-INITVAL}
Proposition \ref{PROPFIR} implies that, in order to prove existence and uniqueness global solutions, 
it suffices to prove the same for the firing times $T_nu$ (see section \ref{S-DEF}). By induction, it suffices in turn to prove existence and uniqueness of the first firing time function $T_1u$ whose equation can be rewritten as follows
\begin{equation}
Mu(x,t)>\eta,\ \forall t\in [0,T_1u(x))\quad\text{and}\quad Mu(x,T_1u(x))=\eta,\ \forall x\in (0,1].
\label{EQFIRING}
\end{equation}
where the quantity $u(x,t)$ is given by \eqref{EXPRSOLUTION} with $n=1$, {\sl viz.}\ 
\[
u(x,t)=\left\{\begin{array}{ccl}
(u(x)-t)^+&\text{if}&0\leq t\leq T_1u(x)\\
(1-t+T_1u(x))^+&\text{if}&T_1u(x)< t
\end{array}\right.\ \forall x\in (0,1],t\in [0,T_1u(1)].
\]
Existence and uniqueness of solutions to \eqref{EQFIRING} 
can be granted in two distinct cases; either when the initial profile is locally constant in a right neighborhood of every point in $(0,1)$, or in the weak coupling regime $\epsilon\leq 1$. 
\begin{Pro}\label{GLOBEXIST}
Let $\eta$ be arbitrary. There exists a unique global solution to equation \eqref{DEFDYNAM} in the following cases
\begin{itemize}
\item[$\bullet$] $\epsilon$ is arbitrary and the profile $u$ is such that there exists $\Delta_x>0$ for every $x\in (0,1)$, so that $u$ is constant on $(x,x+\Delta_x]$,
\item[$\bullet$] $\epsilon\leq 1$ and $u$ is arbitrary. 
\end{itemize}
\end{Pro}
More generally, the proof implies that $T_1u$ can be uniquely determined under the following assumption: {\sl $T_1u$ is either locally constant or it satisfies $T_1u(x)\leq u(x)$, in the right neighborhood of every point}. An example 
is given by periodic trajectories associated with profiles that are strictly increasing in a right neighborhood of every point, see comment after Proposition \ref{PROEXISTPERIOD} and the inequality \eqref{PERIODNODAMP} below. However, there are cases that do not fit this setting and for which proving the existence of global solutions remains open, especially if there exists $x\in [0,1)$ such that 
\[
T_1u(x+0)=u(x+0)\quad\text{and}\quad T_1u(x)<T_1u(y),\ \forall y>x.
\]

\noindent
{\sl Proof of Proposition \ref{GLOBEXIST}.} We consider the two cases separately.

\noindent
$\bullet$ $u$ is locally constant on right neighborhoods. Group invariance (property in section \ref{ONE}) implies that the solution $T_1u$ of equation \eqref{EQFIRING} must also be locally constant on right neighborhoods. In particular, if $u(x)=u(\Delta_0)$ for all $x\in (0,\Delta_0]$, then we must have $T_1u(x)=T_1u(\Delta_0)$ on the same interval. We first aim to determine $T_1u(\Delta_0)$ and more generally, to determine the firing time of the first firing plateau. 

Let the function $S$ be defined by
\[
S(x)= \int_x^1u(y)-u(x)dy=\int_x^1u(y)dy-(1-x)u(x),
\]
and consider separately the two cases $\epsilon S(\Delta_0)\leq 1$ and $\epsilon S(\Delta_0)>1$. 

\noindent
The inequality in the first case is equivalent to $Mu(\Delta_0)-u(\Delta_0)\leq \eta$. Hence, the firing time must be given by  
\[
T_1u(\Delta_0)=Mu(\Delta_0)-\eta\leq u(\Delta_0)\quad \text{and}\quad T_1u(\Delta_0)<T_1u(x),\ \forall x>\Delta_0.
\]
For consistence with the second case, we set $\Delta'_0=\Delta_0$ in this case.

\noindent
In the second case, let 
\[
\Delta'_0:=\sup\{x>0\ :\ \epsilon S(x)>1\}.
\]
Left continuity and boundedness of $u$ imply that $S$ is also left continuous. In addition, we claim that it is non-increasing. (Indeed, using that the derivative of $u$ exists and is finite a.e., it results that the same property holds for $S$ and its derivative is given by $-(1-x)u'(x)\leq 0$.) Hence we have $\Delta'_0\geq \Delta_0$. It is immediate to check that $\Delta'_0$ corresponds to the size of the first firing plateau, {\sl viz.}\ we have 
\[
T_1u(x)=T_1u(\Delta'_0)<u(\Delta'_0+0), \forall x\in (0,\Delta'_0]\quad\text{and}\quad  T_1u(\Delta_0)<T_1u(x),\ \forall x>\Delta'_0.
\]
The number $T_1u(\Delta'_0)$ itself is uniquely defined by 
\[
\epsilon\int_{\Delta'_0}^1u(y)-T_1u(\Delta'_0)dy= 1.
\]

Now, by repeating the same arguments to the translated profile $v(x)$ immediately after $T_1u(\Delta'_0)$ and defined by
\[
v(x)=\left\{\begin{array}{ccl}
u(x+\Delta'_0,T_1u(\Delta'_0)+0)&\text{if}&0< x\leq 1-\Delta'_0\\
1&\text{if}&1-\Delta'_0<x\leq 1
\end{array}\right.
\]
an induction concludes that if $T_1u$ is already defined on $(0,x]$ (where $x$ is arbitrary), then it can be uniquely extended to $(x,x+\Delta'_x]$ where $\Delta'_x\geq \Delta_x$. Moreover, recall that $T_1u$ must be non-decreasing and left continuous. Hence, if it is defined on any set $S$, it can always be uniquely extended to the semi-closed set $S^\ell:=\bigcap_{\delta>0}S+[0,\delta)$. The following technical statement then implies that this process uniquely defines $T_1u$ on $(0,1]$.
\begin{Lem}
Let $S\subset (0,+\infty)$ be a semi-closed set such that 
\begin{itemize}
\item[$\bullet$] there exists $x_0\in (0,1]$ such that $(0,x_0]\subset S$,
\item[$\bullet$] for every $x\in S$, there exists $\Delta_x>0$ such that $(x,x+\Delta_x]\subset S$.
\end{itemize}
Then $(0,1]\subset S$.
\end{Lem}
{\sl Proof of the Lemma.} The proof proceeds by transfinite induction. Starting with $(0,x_0]$, the induction property implies $(0,x_1]\subset S$ where $x_1>x_0+\Delta_{x_0}$, and then 
\[
(0,x_n]\subset S,\ \forall n\in\N.
\]
Let $x_\omega$ be the limit of the increasing sequence $\{x_n\}_{n\in\N}$. By induction, we also have $(0,x_\omega)\subset S$ and then $(0,x_\omega]\subset S^\ell=S$. 

If $x_\omega\geq 1$, we are done. Otherwise, we continue the induction to successive ordinals, until we eventually reach an uncountable ordinal $\omega_1$. Then we must have $x_{\omega_1}>1$, otherwise we would have a uncountable collection of contiguous intervals whose union covers $(0,x_{\omega_1}]$ but not $(0,1]$. This is impossible; hence  
$(0,1]\subset (0,x_{\omega_1}]$ and the proof is complete. \hfill $\Box$
\medskip

\noindent
$\bullet$ $\epsilon\leq 1$. Firing without grouping (property in section \ref{TWOPRIME}) implies $T_1u\leq u$ in this case; hence solution components $u(x,t)$ remain positive for $t\leq T_1u(1)$. In particular, using the definition of the lower trace $\underline{T_1u}$ of the firing profile $T_1u$ (see Appendix \ref{S-TRACE}), 
we get the following expression 
\begin{equation}
u(y,T_1u(x))=\left\{\begin{array}{ccl}
1-T_1u(x)+T_1u(y)&\text{if}&0<y<\underline{T_1u}(x)\\
u(x)-T_1u(x)&\text{if}&\underline{T_1u}(x)=x\ \text{and}\ y=x\\
u(y)-T_1u(x)&\text{if}&\underline{T_1u}(x)<y\leq 1
\end{array}\right.\ \forall x\in (0,1]
\label{EXPRFIRING}
\end{equation}
(There is no need to define $u(\underline{T_1u}(x),T_1u(x))$ when $\underline{T_1u}(x)<x$.) To proceed, we shall need that the profile and firing time traces must coincide in absence of grouping, as now stated.
\begin{Claim}
If $T_1u(x)<T_1u(y)$ for every pair $(x,y)\in (0,1]$ such that $u(x)<u(y)$, then we have $\underline{T_1u}(x)=\underline{u}(x)$ for all $x\in (0,1]$.
\label{EQUTRACES}
\end{Claim}
{\sl Proof of the Claim.} We consider the cases $\underline{u}(x)<x$ and $\underline{u}(x)=x$ separately. 

\noindent
In the first case, the equality $u(y)=u(x)$ for all $y\in (\underline{u}(x),x]$ implies $T_1u(y)=T_1u(x)$ for all $y\in (\underline{u}(x),x]$ and thus $\underline{T_1u}(x)\leq \underline{u}(x)$. By contradiction, if we had $\underline{T_1u}(x)< \underline{u}(x)$, then for any $y\in (\underline{T_1u}(x),\underline{u}(x))$ we would have 
\[
u(y)<u(x)\quad\text{and}\quad T_1u(y)\geq T_1u(x),
\]
which contradicts the assumption of the Claim.

\noindent
For the second case, notice that we have $u(y)<u(x)$ for all $y<\underline{u}(x)$, and then $T_1u(y)<T_1u(x)$ for all $y<\underline{u}(x)$, from the Claim assumption. This implies $\underline{u}(x)\leq \underline{T_1u}(x)$ and the conclusion follows from the facts that $\underline{u}(x)=x$ and $\underline{T_1u}(x)\leq x$. \hfill $\Box$

The inequality $T_1u(x)\leq u(x)$ together with 'firing without grouping' ensures that the assumption in Claim \ref{EQUTRACES} holds. Using also expression \eqref{EXPRFIRING} to manipulate equation \eqref{EQFIRING},
we obtain the following affine functional equation
\begin{equation}
(\text{Id}-L_{\underline{u}})T_1u(x)=(1-\epsilon\eta)u(x)-\eta+\epsilon\eta\left(\underline{u}(x)+\int_{\underline{u}(x)}^1u(y)dy\right),\ \forall x\in (0,1],
\label{IDMINUSL}
\end{equation}
where the linear operator $L_{\underline{u}}$ is defined by
\[
L_{\underline{u}}v(x)=\epsilon\eta\int_0^{\underline{u}(x)}v(y)dy,\ \forall x\in (0,1],
\]
for every bounded Borel measurable function $v$ defined on $(0,1]$. Endowing the corresponding space with the uniform norm $\|\cdot\|_\infty$, the assumption $\epsilon<1/\eta$ implies $\|L_{\underline{u}}\|_\infty<1$. Hence, the operator $\text{Id}-L_{\underline{u}}$ is invertible with bounded inverse. Accordingly, there exists a unique bounded solution $x\mapsto T_1u(x)$ to equation \eqref{IDMINUSL}. \hfill $\Box$

\section{Asymptotic properties of the dynamics}\label{S-ASYMP}
This section investigates the asymptotic behavior of global solutions as $t\to+\infty$. To that goal, we equip the set of bounded Borel measurable functions defined on $(0,1]$, with the $L^1$-norm $\|\cdot\|_1$.

\subsection{Existence condition and uniqueness of periodic trajectories}
Anticipating the results below on asymptotic periodicity, we present here preliminary properties of periodic trajectories. 
By a {\bf periodic trajectory}, we mean a solution $\{u(x,t)\}$ such that there exists $\tau\in \R^+$ so that 
\[ 
u(x,t+\tau+0)=u(x,t),\ \forall x\in (0,1],\ t\in\R^+.
\]
Of note, the period $\tau$ here  
can only be one of the numbers $\{T_nu(1)\}_{n\in\N}$, because the profile $u(\cdot,t)$ cannot be both non decreasing and satisfy $u(1,t+0) = 1$ at other times. Here, we shall focus on $T_1u(1)$-periodic trajectories ({\sl i.e.} $\tau=T_1u(1)$) because these solutions turn out to play a special role in the asymptotic dynamics.

As our next statement indicates, periodic trajectories are uniquely determined by the traces $\underline{u}$ and $\overline{u}$ of their initial profile. Recall from Appendix \ref{S-TRACE} 
that every lower trace function is entirely determined by a countable collection of pairwise disjoint semi-open intervals in $(0,1]$ and the knowledge of a lower trace completely determines the upper trace.
\begin{Pro}
{\rm (i)} Let $\eta,\epsilon$ be any parameters and let $u_\text{tr}$ be any lower trace function. There exists at most one non-decreasing profile $u$ such that $\underline{u}=u_\text{tr}$ and such that the trajectory issued from $u$ is periodic with period $T_1u(1)$.

\noindent
{\rm (ii)} This periodic trajectory exists iff 
\begin{equation}
\epsilon\lesssim \frac{1}{\int_0^1u_\text{tr}(x)dx\left(1-\inf\limits_{x\in (0,u_\text{tr}(1))}\frac{\overline{u_\text{tr}}(x)-u_\text{tr}(x)}{1-\overline{u_\text{tr}}(\overline{u_\text{tr}}(x)+0)+\overline{u_\text{tr}}(x)}\right)},
\label{EXISTPERIOD}
\end{equation}
where the symbol $\lesssim$ means $<$ if the infimum is a minimum, and it means $\leq$ if this bound is not attained. 
\label{PROEXISTPERIOD}
\end{Pro}
The proof is given below. Notice that the condition \eqref{EXISTPERIOD} does not depend on $\eta$ and the infimum is a minimum for every finite step trace function. In more general cases, this infimum may be attained or not. Moreover, the condition simply reduces to $\epsilon\lesssim \frac{1}{\int_0^1u_\text{tr}(x)dx}$ for every trace $u_\text{tr}$ for which for every $\delta>0$, there exists $x_\delta\in (0,1]$ such that 
\[
\overline{u_\text{tr}}(x_\delta)-u_\text{tr}(x_\delta)\leq \delta.
\]
In addition, using $u_\text{tr}(x)\leq x$ and properties of the traces (Appendix \ref{S-TRACE}), 
it is easy to conclude that the denominator in \eqref{EXISTPERIOD} is certainly not larger than $\tfrac12$, and this bound is attained for the increasing trace $u_\text{tr}(x)=x$ for all $x$. Moreover, the quantity in the right hand side of \eqref{EXISTPERIOD} continuously depends on $u_\text{tr}$ ($L^1$-topology). In particular, if all existing plateaus are sufficiently small (depending on $\epsilon$), then the associated periodic orbit does not persist when $\epsilon>2$ is sufficiently large. These facts yield the following statement.
\begin{Cor}
 {\rm (i)} There exists a $T_1u(1)$-periodic trajectory with $\underline{u(\cdot,0)}=u_\text{tr}$ for every lower trace function $u_\text{tr}$, iff $\epsilon<2$.
 
 \noindent
 {\rm (ii)} For every $\epsilon>2$, there exists $\ell_\epsilon>0$ such that, for every lower trace function $u_\text{tr}$ so that $\|\overline{u_\text{tr}}-u_\text{tr}\|_1<\ell_\epsilon$, the associated $T_1u(1)$-periodic trajectory does not exist. 
\label{COREXISTPERIOD}
\end{Cor}
 
\noindent
{\sl Proof of Proposition \ref{PROEXISTPERIOD}.} We study the equation 
\[
u(\cdot,T_1u(1)+0)=u\quad\text{and}\quad \underline{u}=u_\text{tr},
\]
for an arbitrary trace $u_\text{tr}$. Of note, as a profile immediately after firing, the function $u$ must satisfy $u>0$. Equation \eqref{EXPRSOLUTION}, together with this assumption, implies that this equation rewrites as 
\[
1-T_1u(1)+T_1u=u,
\]
which shows that the difference $T_1u-u=T_1u(1)-1$ must be a constant function. This elicits two cases; either this difference is negative, or it is non-negative. 

Assume the first case. Replacing $T_1u$ by $u+T_1u(1)-1$, and $\underline{u}$ by $u_\text{tr}$, in equation \eqref{IDMINUSL} for the first firing time (where now $\Delta_0=1$), one obtains after some simple algebra the following equation for $u$ and $T_1u(1)$:
\[
\epsilon\eta u+\eta=1-T_1u(1)+\epsilon\eta T_1u(1)u_\text{tr}+\epsilon\eta\int_0^1u(y)dy.
\]
Integrating over $(0,1]$ yields the following expression 
\[
T_1u(1)=\frac{1-\eta}{1-\epsilon\eta\int_0^1u_\text{tr}(x)dx}
\]
Moreover, using $u(1)=1$ in the equation for the function $u$ above and evaluated at $x=1$ in order to eliminate the integral term, one gets that the solution writes
\[
u =1-T_1u(1)(u_\text{tr}(1)-u_\text{tr}).
\]
Clearly, this solution is unique and is left continuous. Moreover, it is non-decreasing and we have\footnote{If $u,v$ are two real functions, then $u\leq v$ (resp.\ $u<v$) means $u(x)\leq v(x)$ (resp.\ $u(x)<v(x)$) for all $x\in (0,1]$.} $u\leq 1$ iff $T_1u(1)\geq 0$, {\sl i.e.}\ iff 
\[
\epsilon\eta< \frac{1}{\int_0^1u_\text{tr}(x)dx},
\]
using the constraint $\eta<1$. Finally, we need to make sure that $T_1u-u$ is negative, {\sl viz.}\ $T_1u(1)< 1$ (which also implies $u>0$). This inequality is equivalent to the following one 
\begin{equation}
 \epsilon< \frac{1}{\int_0^1u_\text{tr}(x)dx}.
\label{PERIODNODAMP}
\end{equation}
Altogether, we conclude that a unique initial profile with trace $u_\text{tr}$ generates a $T_1u(1)$-periodic trajectory with $T_1u<u$ iff the inequality \eqref{PERIODNODAMP} holds.

The analysis is similar in the second case. Here, we use the equality $\underline{T_1u}=u_\text{tr}$ and equation \eqref{ALTERTRACE} in Appendix \ref{S-TRACE} 
to obtain $\overline{T_1u}=\overline{u_\text{tr}}$. Accordingly, the solution expression at time $T_1u(x)$ is given by 
\[
u(y,T_1u(x))=\left\{\begin{array}{ccl}
1-T_1u(x)+T_1u(y)&\text{if}&0<y\leq u_\text{tr}(x)\\
0&\text{if}&u_\text{tr}(x)< y\leq \overline{u_\text{tr}}(x)\\
u(y)-T_1u(x)&\text{if}&\overline{u_\text{tr}}(x)<y\leq 1
\end{array}\right.
\]
(The last line obviously does not apply when $\overline{u_\text{tr}}(x)=1$.) Inserting this expression into the equation $Mu(\cdot,T_1u(\cdot))=\eta$, one gets after some simple algebra (using also $u(y)=u(x)$ for all $y\in (\underline{u}(x),\overline{u}(x)]$)
\[
(1-\overline{u_\text{tr}})(1-T_1u(1))+u_\text{tr}-u+\int_0^1u(y)dy=\frac{1}{\epsilon}.
\]
Integrating over $(0,1)$ and using equation \eqref{NORMATRACE} in Appendix \ref{S-TRACE} 
yields in this case
\[
T_1u(1)=\frac{2\int_0^1u_\text{tr}(x)dx-\frac{1}{\epsilon}}{\int_0^1u_\text{tr}(x)dx}.
\] 
As in the first case, by subtracting from the equation its expression for $x=1$, we get that the solution writes
\[
u=(1-\overline{u_\text{tr}})(1-T_1u(1))+1-u_\text{tr}(1)+u_\text{tr}.
\]
The existence of this solution requires $T_1u(1)\geq 1$ ({\sl i.e.}\ $T_1u-u$ is a non-negative function); this inequality is equivalent to 
\[
\epsilon\geq \frac{1}{\int_0^1u_\text{tr}(x)dx}
\]
which is complementary to the existence condition \eqref{PERIODNODAMP} in the first case. We also need to impose $u>0$. By monotonicity, this condition is equivalent to 
\[
T_1u(1)-1\lesssim \frac{1-u_\text{tr}(1)}{1-\overline{u_\text{tr}}(0+0)}.
\]
Finally, we need to make sure that $\underline{T_1u}=u_\text{tr}$, {\sl viz.}\ we must have 
\[
T_1u(x)<u(y),\ \forall y>\overline{u_\text{tr}}(x), x\in (0,u_\text{tr}(1)).
\]
Using $T_1u(x)=u(x)+T_1u(1)-1$, the expression of $u$ and relation \eqref{RIGHTLIMIT} in Appendix \ref{S-TRACE}, 
this condition is equivalent to 
\[
T_1u(1)-1\lesssim \inf_{x\in (0,u_\text{tr}(1))}\frac{\overline{u_\text{tr}}(x)-u_\text{tr}(x)}{1-\overline{u_\text{tr}}(\overline{u_\text{tr}}(x)+0)+\overline{u_\text{tr}}(x)}.
\]
We have $\inf\limits_x\overline{u_\text{tr}}(x)-u_\text{tr}(x)\leq 1-u_\text{tr}(1)$ and $\inf\limits_x \overline{u_\text{tr}}(\overline{u_\text{tr}}(x)+0)-\overline{u_\text{tr}}(x)\leq \overline{u_\text{tr}}(0+0)$; hence the last condition is the strongest one. Simple algebra finally yields the inequality \eqref{EXISTPERIOD} and this completes the proof. \hfill $\Box$

\subsection{Asymptotic periodicity for solutions with finite step profiles}
Proposition \ref{GLOBEXIST} implies that a global solution of equation \eqref{DEFDYNAM} exists for every initial finite step profile (because every finite step function is obviously locally constant in the right neighborhood of every point). Moreover, the grouping properties in section \ref{S-CLUSTER} imply that the number of profile steps either remains constant or decreases before each firing; hence this number must eventually become constant. Therefore, when dealing with asymptotic properties, we may assume without loss of generality that this number remains constant, {\sl viz.}\ that the lower trace function is periodic after each full cycle of firing.
\begin{Thm}
Let $\eta,\epsilon$ be any parameters and let $u$ be any initial finite step function for which the number of clusters remains constant in the corresponding trajectory. Then, we have
\[
\lim_{t\to+\infty}\|u(\cdot,t+T_1u_\text{per}(1)+0)-u(\cdot,t)\|_1=0,
\]
where $u_\text{per}$ is the periodic trajectory profile such that $\underline{u_\text{per}}=\underline{u}$. 
If, in addition, the condition \eqref{EXISTPERIOD} holds with $u_\text{tr}=\underline{u}$, then we also have 
\[
\lim_{n\to+\infty}\|u(\cdot,T_n u(1)+0)-u_\text{per}\|_1=0.
\]
\label{GLOBCONV}
\end{Thm}
It may happen that the first limit in this statement holds although the condition \eqref{EXISTPERIOD} fails, more precisely when 
\[
\epsilon= \frac{1}{\int_0^1u_\text{tr}(x)dx\left(1-\min\limits_{x\in (0,u_\text{tr}(1))}\frac{\overline{u_\text{tr}}(x)-u_\text{tr}(x)}{1-\overline{u_\text{tr}}(\overline{u_\text{tr}}(x)+0)+\overline{u_\text{tr}}(x)}\right)}.
\]
In this case, the periodic cycle $\{u_\text{per}(\cdot,t)\}$ is an example of ghost orbit mentioned in Appendix \ref{S-DISCONT}. Besides, the proof of the Theorem below implies that when
\[
\epsilon> \frac{1}{\int_0^1u_\text{tr}(x)dx\left(1-\min\limits_{x\in (0,u_\text{tr}(1))}\frac{\overline{u_\text{tr}}(x)-u_\text{tr}(x)}{1-\overline{u_\text{tr}}(\overline{u_\text{tr}}(x)+0)+\overline{u_\text{tr}}(x)}\right)},
\]
any solution issued from a finite step profile $u$ such that $\underline{u}=u_\text{tr}$, must experience the grouping of a least two plateaus in finite time. In this way, we have obtained a complete picture of all possible asymptotic finite-dimensional population distributions, depending on the coupling intensity $\epsilon$.  
\smallskip

\noindent
{\sl Proof of Theorem \ref{GLOBCONV}.} The dynamics of finite step functions is entirely determined by the number $N$ of steps, by the step lengths $\{\ell_n\}_{n=1}^N$ and by the step expression levels $\{u_n\}_{n=1}^N$. (Here, step labelling follows from cell labelling on $(0,1]$, {\sl viz.}\ $n=1$ means the first group $x\in (0,x_1]$ for some $x_1>0$, $n=2$ means $x\in (x_1,x_2]$ with $x_2>x_1$ and so on.) 

Assuming that no grouping occurs in time implies that not only $N$ but also the length distribution $\{\ell_n\}_{n=1}^N$ remains constant. Only the expression levels $\{u_n\}_{n=1}^N$ may depend on time. In this setting, we aim to prove that the time evolution of every expression level vector is asymptotically periodic in $\R^N$. 

To that goal, it suffices to consider the discrete time dynamics that brings the expression levels after a firing to those after the next firing (Poincar\'e return map). From thereon in this proof, by {\bf time}, we mean the integer $t\in\N$ that labels the vector $\{u_n^t\}_{n=1}^N$ after the $t^\text{th}$ firing.

It turns out more convenient to combine the discrete time dynamics with the cyclic permutation of indices, so that any vector with non-decreasing components and $u_N=1$ is mapped after every iteration onto a vector carrying the same properties. Ignoring systematically the last component $u_N=1$, this amounts to considering iterations of the $(N-1)$-dimensional map implicitly defined by $\{u_n^t\}_{n=1}^{N-1}\mapsto \{u_{n+1}^{t+1}\}_{n=1}^{N-1}$. 

Beside the original parameters $\eta$ and $\epsilon$, this firing map $F_{\ell}$ is also parametrized by the step length distribution $\ell:=\{\ell_n\}_{n=1}^{N}$ (at time $t$). Its action on vectors $u:=\{u_n\}_{n=1}^{N-1}\in {\cal U}_{N-1}$, where 
\[
{\cal U}_{N-1}=\left\{u=\{u_n\}_{n=1}^{N-1}\ :\ 0< u_1< u_2< \cdots <u_{N-1}<1\right\},
\]
is explicitly given by 
\[
(F_\ell u)_n=u_{n+1}-T_\ell u,\ \forall n=1,\cdots, N-1,
\]
where $T_\ell u$ is the time of the first firing in the trajectory starting from the finite step profile associated with $u$ (at time $t$). Iterations of the firing map have to incorporate permutations of the step length distribution, {\sl i.e.}\ we need to consider the composed map\footnote{If the length distribution period $N_\text{per}$, defined by 
\[
N_\text{per}=\min\left\{k\ :\ \ell_{n+k\ \text{mod}\ N}=\ell_{n},\ \forall n\in \{1,\cdots,N\}\right\},
\]
happens to be smaller than $N$, it actually suffices to consider the composed map $F_{R^{N_\text{per}-1}\ell}\circ \cdots \circ F_{R\ell}\circ F_{\ell}$, because appropriate permutations of the profiles associated to iterates of this map are non-increasing finite step functions with identical step length distribution as the initial profile. The same consideration suggests to consider, given any step length distribution, the permutation that minimises the period $N_\text{per}$.} 
\begin{equation}
F^N_\ell:=F_{R^{N-1}\ell}\circ  F_{R^{N-2}\ell}\circ \cdots \circ F_{R\ell}\circ F_{\ell},
\label{COMPOSED}
\end{equation}
where $(R\ell)_n=\ell_{n+1\ \text{mod}\ N}$ for all $n\in \{1,\cdots ,N\}$. By identifying each $\{u_n\}_{n=1}^{N-1}\in {\cal U}_{N-1}$ with its corresponding $N$-step profile, the map $F^N_\ell$ is clearly equivalent to the return map $u(\cdot,0)\mapsto u(\cdot,T_1u(1)+0)$. Theorem \ref{GLOBCONV} is therefore an immediate consequence of the following statement. Let $\|\cdot\|$ be any norm in $\R^{N-1}$.
\begin{Pro}
For every step length distribution $\ell=\{\ell_n\}_{n=1}^N$, there exist $C>0$ and $\rho\in [0,1)$ such that for every pair $u,v\in {\cal U}_{N-1}$ for which $(F_\ell^N)^ku,(F_\ell^N)^kv\in {\cal U}_{N-1}$ for all $k\in\N$, we have 
\begin{equation}
\|(F_\ell^N)^ku-(F_\ell^N)^kv\|\leq C\rho^k\|u-v\|,\ \forall k\in\N.
\label{EXPONBOUND}
\end{equation}
\label{PROBOUND}
\end{Pro}
{\sl Proof of the Proposition.} 
Independently of considerations on parameters ({\sl i.e.}\ whether or not $\epsilon\leq 1$), there are {\sl a priori} two cases for the expression of $T_\ell u$, depending on whether $u_1$ fires before reaching 0 or after. Simple calculations yield the following expression: 
\[
T_\ell u=\left\{\begin{array}{ccl}
(1-\epsilon\eta)u_1+\epsilon\eta\sum\limits_{n=1}^N\ell_nu_n-\eta&\text{if}&u_1\geq T_\ell u\\
\frac{1}{1-\ell_1}\left(\sum\limits_{n=2}^N\ell_nu_n-\frac{1}{\epsilon}\right)&\text{if}&u_1\leq T_\ell u
\end{array}\right.
\]
Accordingly, for every $\ell$, the map $F_{\ell}$ is a continuous (piecewise) affine map with at most two pieces, say $F_\ell^+$ and $F_\ell^-$. Let $DF_\ell^+$ and $DF_\ell^-$ be the linear maps associated with each piece. We claim that, for every vector pair $u,v\in {\cal U}_{N-1}$, there exists $\alpha\in [0,1]$ such that 
\begin{equation}
F_\ell u-F_\ell v=(1-\alpha)DF_\ell^+(u-v)+\alpha DF_\ell^-(u-v).
\label{CONVEXSUM}
\end{equation}
Equation \eqref{CONVEXSUM} is obvious (and holds with $\alpha\in \{0,1\}$) if there exists $s\in \{+,-\}$ such that $F_\ell w=F_\ell^s w$ for $w=u,v$. Otherwise, by continuity of the firing time function $w\mapsto T_\ell w$, there exists $\beta\in (0,1)$ such that 
\[
F_\ell (\beta u+ (1-\beta)v)=F_\ell^+ (\beta u+ (1-\beta)v)=F_\ell^-(\beta u+ (1-\beta)v).
\]
If $F_\ell u=F_\ell^+ u$ and $F_\ell v=F_\ell^- v$, one gets 
\begin{align*}
F_\ell u-F_\ell v&=F_\ell^+u -F_\ell^+ (\beta u+ (1-\beta)v)+F_\ell^- (\beta u+ (1-\beta)v)-F_\ell^- v\\
&=(1-\beta)DF_\ell^+(u-v)+\beta DF_\ell^-(u-v)
\end{align*}
{\sl viz.}\ equation \eqref{CONVEXSUM} holds with $\alpha=\beta$. Otherwise, we must have $F_\ell u=F_\ell^- u$ and $F_\ell v=F_\ell^+ v$, and a similar calculation shows that equation \eqref{CONVEXSUM} holds with $\alpha=1-\beta$.

The assumption $(F_\ell^N)^ku,(F_\ell^N)^kv\in {\cal U}_{N-1}$ for all $k\in\N$ implies that we must have 
\[
F_{R^{j}\ell}\circ \cdots \circ F_{\ell}\circ (F_\ell^N)^ku\quad\text{and}\quad F_{R^{j}\ell}\circ \cdots \circ F_{\ell}\circ(F_\ell^N)^kv\in {\cal U}_{N-1},\ \forall j\in\{0,\cdots , N-1\},\ k\in\N.
\]
By induction, the composed map $F_\ell^N$ defined by \eqref{COMPOSED} is also a continuous (piecewise) affine map with {\sl a priori} $2^{N-1}$ pieces; each piece writes
\[
F_{R^{N-1}\ell}^{s_{N-1}}\circ  F_{R^{N-2}\ell}^{s_{N-2}}\circ \cdots \circ F_{R\ell}^{s_1}\circ F_{\ell}^{s_0}
\]
for some $\{s_n\}_{n=0}^{N-1}\in \{-,+\}^N$. Moreover, equation \eqref{CONVEXSUM} implies that, for every vector pair $u,v$, there exists $\alpha_{\{s_n\}_{n=0}^{N-1}}\in [0,1]$ for every piece, such that $\sum\limits_{\{s_n\}_{n=0}^{N-1}\in \{-,+\}^N}\alpha_{\{s_n\}_{n=0}^{N-1}}=1$ and 
\begin{equation}
F_\ell^Nu-F_\ell^Nv=\sum\limits_{\{s_n\}_{n=0}^{N-1}\in \{-,+\}^N}\alpha_{\{s_n\}_{n=0}^{N-1}} DF_{R^{N-1}\ell}^{s_{N-1}} \cdots DF_{\ell}^{s_0}(u-v).
\label{DIFFITER}
\end{equation}

Consider now the joint spectral radius (see {\sl e.g.}\ \cite{RS60}) of the finite collection 
\[
{\cal A}=\{DF_{R^{j}\ell}^{s}\}_{j\in \{0,\cdots,N-1\},s\in\{-,+\}}
\]
defined by
\[
\text{JoinSpecRad}({\cal A}):=\limsup_{k\to+\infty}\left(\max\left\{\left\|\prod_{i=1}^kA_i\right\|\ :\ A_i\in {\cal A},\ \forall i\in\{1,\cdots ,k\}\right\}\right)^\frac{1}{k}.
\]
By applying relation \eqref{DIFFITER} repeatedly to the iterates $(F_\ell^N)^ku$ and $(F_\ell^N)^kv$, a straightforward argument shows that the inequality \eqref{EXPONBOUND} holds under the condition $\text{JoinSpecRad}({\cal A})<1$. Furthermore, since ${\cal A}$ is a finite set, we have \cite{BW92}
\[
\text{JoinSpecRad}({\cal A})=\limsup_{k\to+\infty}\left(\max\left\{\text{SpecRad}\left(\prod_{i=1}^kA_i\right)\ :\ A_i\in {\cal A},\ \forall i\in\{1,\cdots ,k\}\right\}\right)^\frac{1}{k},
\]
where $\text{SpecRad}$ is the standard spectral radius of a matrix. Therefore, in order to prove the Proposition, it suffices to show  that the right hand side here is smaller than 1.  

To proceed, we observe that, by applying the change of variable $\{u_n\}_{n=1}^{N-1}\mapsto \{v_n\}_{n=1}^{N-1}$, where 
\[
v_n=\left\{\begin{array}{ccl}
u_n-u_{n+1}&\text{if}&n\in \{1,\cdots,N-2\}\\
u_{N-1}&\text{if}&n=N-1
\end{array}\right.,
\]
the matrix $DF_{\ell}^{s}$ transforms into $A_{a(\ell,s)}$, where the vectors $a(\ell,s)=\{a_n(\ell,s)\}_{n=1}^{N-1}$ for $s=-,+$ are respectively defined by 
\[
a_n(\ell,+)=1-\epsilon\eta\sum_{m=n+1}^{N}\ell_m\quad \text{and}\quad a_n(\ell,-)=\frac{1}{1-\ell_1}\sum\limits_{m=2}^{n+1}\ell_m,\  \forall n=1,\cdots N-1,
\]
and where, given an arbitrary vector $a=\{a_k\}_{k=1}^K$ ($K\in\N$), the matrix $A_{a}$ is the following $K\times K$ companion matrix
\[
A_{a}=\left(\begin{array}{ccccc}
0        &1        &0         &\cdots&0\\
\vdots&\ddots&\ddots&\ddots&\vdots\\
\vdots &          &\ddots&\ddots&0\\
0         &\cdots&\cdots&0&1\\
-a_{1}&-a_{2}&\cdots&\cdots&-a_{K}
\end{array}\right).
\]
Using that the coefficients $a_n(\ell,+)$ and $a_n(\ell,-)$ are positive and decreasing, a classical result in numerical analysis \cite{J64} implies that $\text{SpecRad}(A_{a(\ell,s)})<1$ for $s\in\{-,+\}$, for every length distribution $\ell=\{\ell_n\}_{n=1}^{N-1}$. However, the products $\prod\limits_{i=1}^kA_{a(\ell_i,s_i)}$ of such matrices usually do not commute (and the associated semi-group is usually infinite) and do not appear to have any special form that would allow one 
to immediately conclude about their spectral radius. Instead, we shall need the following statement that takes advantage of the matrix characteristics.
\begin{Lem} {\rm \cite{KV12}}
Let $\{a^{(j)}\}_{j=1}^{J}$ be a finite collection of $K$-dimensional vectors $a^{(j)}=\{a_k^{(j)}\}_{k=1}^K$ whose components satisfy the following constraint
\[
\min\limits_{k\in\{1,\cdots ,K\}}a_k^{(j)}>0,\ \forall j\in \{1,\cdots ,J\}.
\]
Then, for any collection $\{j_i\}_{i=1}^k$ with $j_i\in \{1,\cdots ,J\}$ for all $i$, we have 
\[
\text{\rm SpecRad}\left(\prod_{i=1}^kA_{a^{(j_i)}}\right)\leq \left(\max\limits_{j\in \{1,\cdots,J\}}\max\limits_{k\in \{1,\cdots,K\}}\frac{a_k^{(j)}}{a_{k+1}^{(j)}}\right)^k
\]  
where we have set $a_{K+1}^{(j)}:=1$ for all $j$. 
\end{Lem}

Applying the Lemma to the collection $\{a(R^j\ell,s)\}_{j\in\{0,\cdots,N-1\},s\in\{-,+\}}$ and letting $a_{N}(R^j\ell,s)=1$, we immediately conclude
\[
\text{JoinSpecRad}({\cal A})=\max_{j\in\{0,\cdots,N-1\},s\in\{-,+\},n\in\{1,\cdots,N-1\}}\frac{a_n(R^j\ell,s)}{a_{n+1}(R^j\ell,s)}<1,
\]
as required. \hfill $\Box$

\subsection{Asymptotic periodicity in the weak coupling regime}
In this section, we go back to arbitrary non-decreasing left continuous initial profiles. Let $\mu_c\in (0.46,0.47)$ be the positive solution of the equation
\[
e^\mu+\frac{\mu^2}{1-\mu}=2.
\]
\begin{Pro}
Let $\eta$ and $\epsilon$ be such that $\epsilon\eta<\mu_c$. There exists $\rho\in [0,1)$ so that for every pair of initial profiles $u,v$ satisfying $\underline{u}=\underline{v}$ and
\[
T_1u\leq u\quad \text{and}\quad T_1u(1)<T_2u\quad (\text{resp.}\quad T_1v\leq v\quad \text{and}\ T_1v(1)<T_2v),
\]
we have
\[
\left\|u(\cdot,T_1u(1)+0)-v(\cdot,T_1v(1)+0)\right\|_1\leq \rho\left\|u-v\right\|_1.
\]
\label{CONVWEAK}
\end{Pro}
By choosing $v=u(\cdot,T_1u(1)+0)$, the previous condition obviously holds for every trajectory for which no grouping ever happens. As the next statement says, this is the case of every trajectory in the weak coupling regime and the conclusion holds provided that $\eta<\mu_c$ (recall that this threshold is assumed to be small from the experiment under modeling). 
\begin{Cor}
Let $\eta<\mu_c$ and $\epsilon<1$. For every initial profile, the subsequent trajectory is asymptotically periodic. 
\end{Cor}

More generally, Proposition \ref{CONVWEAK} implies asymptotic stability of periodic trajectories (associated with infinite dimensional profiles). Indeed, the periodic trajectory associated with the infinitely many step profile with trace $u_\text{tr}$ is known to exist iff
\[
\epsilon\int_0^1u_\text{tr}\leq 1.
\]
Moreover, the analysis in the proof of Proposition \ref{PROEXISTPERIOD} shows this condition implies the one in Proposition \ref{CONVWEAK}. Accordingly for $\epsilon\eta<\mu_c$, this trajectory attracts all solutions in its neighborhood whose profiles after every cycle of firing have trace $u_\text{tr}$. 
\smallskip

\noindent
{\sl Proof of the Proposition.} The assumption $T_1u(1)<T_2u$ implies that the solution $u(\cdot,T_1u(1)+0)$ after a full cycle of firing is given by
\[
u(x,T_1u(1)+0)=1-T_1u(1)+T_1u(x),\ \forall x\in (0,1].
\]
Similarly, we have $v(x,T_1v(1)+0)=1-T_1v(1)+T_1v(x)$ for all $x$. The assumption $T_1u\leq u$ implies that the first firing time $T_1u$ is defined by expression \eqref{IDMINUSL} in Section \ref{S-INITVAL}. A similar expression holds for $T_1v$. In addition, the assumption $\underline{u}=\underline{v}$ implies $L_{\underline{u}}=L_{\underline{v}}$. 

\noindent
By inverting the operator $\text{Id}-L_{\underline{u}}$ by means of a Neumann series, and letting $\mu=\epsilon\eta$, we obtain after some simple algebra
\begin{align*}
u(x,T_1u(1)+0)=&(1-\mu)u(x)\\
&+\mu\sum_{k=1}^{+\infty} L_{\underline{u}}^k\left(\int_0^1u(y)dy-u\right)(x)-L_{\underline{u}}^k\left(\int_0^1u(y)dy-u\right)(1)\\
&+C(x),
\end{align*}
where $C(x)$ does not depend on $u$. A similar expression holds for $v(x,T_1v(1)+0)$. 

To obtain an upper bound on $\left\|u(\cdot,T_1u(1)+0)-v(\cdot,T_1v(1)+0)\right\|_1$, we need to estimate the norm $\|L_{\underline{u}}^kw-L_{\underline{u}}^kw(1)\|_1$ for every measurable function $w$ defined on $(0,1]$ and every $k\in\N$. First, notice that for such functions, by reversing the order of integration, we obtain
\[
\|L_{\underline{u}}w\|_1
\leq\mu\int_0^1\int_0^{\underline{u}(x)}\left|w(y)\right|dydx
\leq \mu\int_0^{\underline{u}(1)}\left(\int_{\overline{u}(y)}^{1}dx\right)|w(y)|dy\leq \mu\|w\|_1.
\]
Moreover, a similar argument implies the following inequality
\[
|L_{\underline{u}}^{k+1}w(x)-L_{\underline{u}}^{k+1}w(1)|\leq \mu\|L_{\underline{u}}^kw\|_1,\ \forall x\in (0,1], k\in\N,
\]
and then the estimate $\|L_{\underline{u}}^kw-L_{\underline{u}}^kw(1)\|_1\leq \mu^{k+1}\|w\|_1$ for all $k\in\N$, by using induction. For constant functions $w(x)=w$, this estimate can be improved by using a better estimate on $\|L_{\underline{u}}w\|_1$. Indeed, using $\underline{u}(x)\leq x$ in a induction implies
\[
\|L_{\underline{u}}^kw\|_1\leq \frac{\mu^k}{(k+1)!}\|w\|_1.
\]
Altogether, we obtain the inequality 
\[
\left\|u(\cdot,T_1u(1)+0)-v(\cdot,T_1v(1)+0)\right\|_1\leq \left(1-\mu+\mu(e^\mu-\mu-1)+\frac{\mu^2}{1-\mu}\right)\|u-v\|_1,
\]
and the Proposition follows from the fact that $1-\mu+\mu(e^\mu-\mu-1)+\frac{\mu^2}{1-\mu}<1$ when $\mu<\mu_c$. \hfill $\Box$

\section{Concluding remarks}
The model equation \eqref{DEFDYNAM} has been inspired by the delayed-differential equation in \cite{MBHT09} that mimics experimental synchronized oscillations reported in \cite{DM-PTH10}. Its features have been selected in order to combine biological relevance, mathematical rigour and non-trivial phenomenology depending on parameters. While some of these features may appear specific, {\sl e.g.}\ permutation symmetry of the interaction term in the repressor definition, possible temporary vanishing of expression levels, instantaneous firings, not to mention the absence of heterogeneities or noise fluctuations, we would like to stress that the phenomenology here does not depend on these features and can be extended by continuity.

Indeed, the contraction property of (sufficiently high iterates of) the return maps after full firing cycles is robust to smooth perturbations of the dynamics. It means for instance that any sufficiently small $C^1$ perturbation of the map $F^N_\ell$, where $\ell$ is an arbitrary finite dimensional step length distribution, satisfies the contraction property \eqref{EXPONBOUND}. Hence, if the associated periodic trajectory exists ({\sl i.e.}\ condition \eqref{EXISTPERIOD} holds for the increasing profile associated with $\ell$), then this trajectory can be continued as a stable periodic orbits of the $C^1$-perturbed map. As a consequence, out of bifurcation points, all attractive periodic trajectories of finite populations persist for sufficiently small perturbation of the original dynamics, independently of the mentioned specific properties of equation \eqref{DEFDYNAM}.

Naturally, this does not precludes a separate investigation into more general degrade-and-fire dynamics. This will be the subject of future studies.
\bigskip
 
\noindent
{\bf Acknowledgments}

\noindent
We are grateful to E.\ Key and H.\ Volkmer for promptly providing us a proof of Lemma 5.5.

\appendix

\section{Discontinuous dependence on initial conditions}\label{S-DISCONT}
Instantaneous resetting simplifies the analysis of equation \eqref{DEFDYNAM}. However, it makes the proof of global existence of solutions rather delicate (and apparently unaccessible by standard approaches such as the Picard operator). It also implies that the solution dependence on initial profiles has discontinuities. In particular, this is the case for the first firing time function $T_1u$.

\noindent
Indeed, there are examples of sequences $\{u_n\}$ of profiles that uniformly converge to a limit profile $u_\infty$ and for which we have $\lim\limits_{n\to+\infty}T_1u_n(x)\neq T_1u_\infty(x)$ for some $x\in (0,1]$.

\noindent
To see this, let $\epsilon\leq 1$, let $u_\infty$ be a profile with a left plateau, {\sl i.e.}\ $u_\infty(x)=u_\infty(x_1)$ for all $x\in (0,x_1]$ ($x_1>0$), and let an approximating sequence be defined by 
\[
u_n(x)=\left\{\begin{array}{ccl}
u_\infty(x)-\frac{1}{n}&\text{if}&0<x\leq \frac{x_1}{2}\\
u_\infty(x)&\text{if}&\frac{x_1}{2}<x\leq 1
\end{array}\right.
\]
We obviously have $T_1u_\infty(x)=T_1u_\infty(x_1)$ for all $x\in (0,x_1]$ and direct calculations yield the following result
\[
\lim_{n\to+\infty}T_1u_n(x)=\left\{\begin{array}{l}
T_1u_\infty(x),\ \forall x\in (0,\frac{x_1}{2}]\\
T_1u_\infty(x)+\epsilon\eta\frac{x_1}{2}(T_1u_\infty(x_1)-u_\infty(x_1)+1),\ \ \forall x\in (\frac{x_1}{2},x_1]
\end{array}\right.
\]
and the inequalities $0\leq u_\infty(x_1)-T_1u_\infty(x_1)<1$ imply $\lim\limits_{n\to+\infty}T_1u_n(x)> T_1u_\infty(x)$ for all $x\in  (\frac{x_1}{2},x_1]$. 

In addition, discontinuities may also result in the existence of attracting {\bf ghost orbits}, depending on parameters. Ghost orbits are periodic cycles of profiles, {\sl viz.}\ $\{u(x,t)\}_{(x,t)\in (0,1]\times \R^+}$ with $u(\cdot,t+\tau+0)=u(\cdot,t)$ for some $\tau>0$, which, while they do not satisfy equation \eqref{DEFDYNAM}, attract all trajectories in their neighborhood (uniform topology). As shown after Theorem \ref{GLOBCONV}, ghost orbits exist at bifurcation points in the parameter space, when a periodic orbit collapses.  

\section{The lower and upper traces of a non-decreasing function}\label{S-TRACE}
Let $u:(0,1]\to (0,1]$ be a left continuous and non-decreasing function. Its lower trace $\underline{u}$ and respectively upper trace $\overline{u}$ are defined as follows  
\[
\underline{u}(x)=\inf\left\{y\in (0,1]\ :\ u(y)\geq u(x)\right\},\ \forall x\in (0,1],
\]
and
\[
\overline{u}(x)=\sup\left\{y\in (0,1]\ :\ u(y)\leq u(x)\right\},\ \forall x\in (0,1].
\]
These functions satisfy the following basic properties. 
\begin{itemize}
\item[$\bullet$] $0\leq \underline{u}(x)\leq x\leq \overline{u}(x)\leq 1$ for all $x\in (0,1]$.
\item[$\bullet$] either $\underline{u}(x)< x$ or $x<\overline{u}(x)$ implies $u(y)=u(x)$ for all $y\in (\underline{u}(x),\overline{u}(x)]$.
\item[$\bullet$] If $u$ is strictly increasing, then $\underline{u}(x)=\overline{u}(x)=x$ for all $x\in (0,1]$.\footnote{Notice that the lower trace can be alternatively defined as $u_\text{inf}^{-1}\circ u$ where the generalized inverse (also called the quantile function in Probability Theory) $u_\text{inf}^{-1}$ can be defined as
\[
u_\text{inf}^{-1}=\inf\left\{y\in (0,1]\ :\ u(y)\geq x\right\},\ \forall x\in (0,1].
\]
In this viewpoint, the property in this item reads $u_\text{inf}^{-1}\circ u=\text{Id}$ for every strictly increasing function $u$. A similar comment applies to the upper trace.}
\item[$\bullet$] $\underline{u}\circ \underline{u}(x)=\underline{u}(x)$ iff $u$ is continuous at $x$.
\item[$\bullet$] $\overline{u}\circ \overline{u}=\overline{u}$.
\item[$\bullet$] Both functions $\underline{u}$ and $\overline{u}$ are left continuous and non-decreasing. (We prove the property for  $\underline{u}$ here; the proof for $\overline{u}$ is similar and is left to the reader. Monotonicity is obvious and implies $\underline{u}(x-0)\leq \underline{u}(x)$. Left continuity is also evident in the case $\underline{u}(x)<x$. If, otherwise $\underline{u}(x)=x$, there must be a sequence $\{x_n\}_{n\in\N}$ such that $u(x_n)<u(x_{n+1})$ and $\lim\limits_{n\to+\infty} x_n=x$. The former condition implies $\underline{u}(x_n)>x_{n-1}$. Together with the latter, we obtain $\underline{u}(x-0)\geq x=\underline{u}(x)$ as desired.)
\end{itemize}
In our context, the traces provide information about the group structure of a population at time $t$: $\underline{u}(x,t)=\overline{u}(x,t)=x$ means that cell $x$ is isolated, while $\underline{u}(x,t)<\overline{u}(x,t)$ means that all cells $y\in (\underline{u}(x,t),\overline{u}(x,t)]$ belong to the same group.

The properties of the lower trace above imply that this function can be entirely determined by its plateaus; namely by considering the following decomposition 
\[
(0,1]={\cal C}_<\cup {\cal C}_=,
\]
where 
\[
{\cal C}_<=\left\{x\in (0,1]\ : \underline{u}(x)<x\right\}
\quad\text{and}\quad
{\cal C}_==\left\{x\in (0,1]\ : \underline{u}(x)=x\right\},
\]
the second item above imposes the existence of a countable (possibly empty) set ${\cal D}$ such that ${\cal C}_<=\bigcup\limits_{i\in {\cal D}}(x_i^{-},x_i^+]$ where $x_i^-<x_i^+\leq x_{i+1}^-$ for all $i$. (Notice that ${\cal C}_=$ is empty when $u$ (or $\underline{u}$) is a step function.) In other words, every countable (possibly empty) collection of pairwise disjoint semi-open intervals in $(0,1]$ uniquely defines a lower trace function. 

The upper trace function depends only on the lower trace, {\sl i.e.}\ $\overline{u}=\overline{u}\circ \underline{u}$ (and {\sl vice-versa}, we have $\underline{u}=\underline{u}\circ \overline{u}$). One can prove this fact using the sets ${\cal C}_<$ and ${\cal C}_=$ and the analogous decomposition for the upper trace. However, for our purpose, it is more convenient to use the following characterization 
\begin{equation}
\overline{u}(x)=\inf\left\{\underline{u}(y)\ :\ x<\underline{u}(y)\right\},\ \forall x\in (0,1],
\label{ALTERTRACE}
\end{equation}
with the convention that $\inf\emptyset=1$ in this expression. To prove this relation, notice first that we must have $\overline{u}(x)\leq\inf\left\{\underline{u}(y)\ :\ x<\underline{u}(y)\right\}$. Indeed, otherwise there existed $y$ such that $x<\underline{u}(y)$ and $\underline{u}(y)<\overline{u}(x)$. Using that the former inequality is equivalent to $\overline{u}(x)<y$, it results that we must have
\begin{equation}
\underline{u}(y)<\overline{u}(x)<y,
\label{IMPOSSIBL}
\end{equation}
which is clearly incompatible with the definition of the traces. Secondly, still by using contradiction, assume that $\overline{u}(x)<\inf\left\{\underline{u}(y)\ :\ x<\underline{u}(y)\right\}$. This implies the existence of $z$ such that $\overline{u}(x)<z<\inf\left\{\underline{u}(y)\ :\ x<\underline{u}(y)\right\}$. However, the first inequality here implies $u(x)<u(z)$ and then $x<\underline{u}(z)$ which contradicts the second inequality. 

Similar arguments prove the following relation 
\begin{equation}
\overline{u}(x)=\underline{u}(\overline{u}(x)+0)=\inf\{\underline{u}(y)\ :\ \overline{u}(x)<y\},\ \forall x\in (0,\underline{u}(1)).
\label{RIGHTLIMIT}
\end{equation}
Indeed, as before, we must have $\overline{u}(x)\leq \inf\{\underline{u}(y)\ :\ \overline{u}(x)<y\}$ because the converse would yield to the double inequality \eqref{IMPOSSIBL} otherwise. Now if there existed $z$ such that 
\[
\overline{u}(x)<z< \inf\{\underline{u}(y)\ :\ \overline{u}(x)<y\},
\]
the right inequality would imply $z<\underline{u}(y)\leq y$ for all $\overline{u}(x)<y$, hence $z\leq \overline{u}(x)$ holds, which contradicts the left inequality.

In the main text, we also refer to the following relation
\begin{equation}
\int_0^1\underline{u}(x)dx+\int_0^1\overline{u}(x)dx=1.
\label{NORMATRACE}
\end{equation}
In order to prove this relation, consider again the decomposition ${\cal C}_<\cup {\cal C}_=$ with ${\cal C}_<=\bigcup\limits_{i\in {\cal D}}(x_i^{-},x_i^+]$. A moment's reflexion yields 
\[
\int_0^1\underline{u}(x)dx=\sum_{i\in {\cal D}}x_i^-(x_i^+-x_i^-)+\int_{\cal C_=}xdx\quad \text{and}\quad
\int_0^1\overline{u}(x)dx=\sum_{i\in {\cal D}}x_i^+(x_i^+-x_i^-)+\int_{\cal C_=}xdx,
\]
and then
\[
\int_0^1\underline{u}(x)dx+\int_0^1\overline{u}(x)dx=\sum_{i\in {\cal D}}(x_i^+)^2-(x_i^-)^2+2\int_{\cal C_=}xdx.
\]
Equation \eqref{NORMATRACE} then directly follows from the fact that  
\[
(x_i^+)^2-(x_i^-)^2=2\int_{x_i^{-}}^{x_i^+}xdx.
\]

\end{document}